\documentclass{ws-procs9x6}
\usepackage{hyperref}
\usepackage{amsfonts}
\input{cyr.fd}
\newtheorem{thm}{Theorem}
\newtheorem{cor}
{Corrolary}

\newtheorem{defn}
{Definition}
\newtheorem{rem}
{Remark}
\providecommand{\myeprint}[2]{E-print: \href{#1}{\texttt{#2}}}
\providecommand{\Space}[3][]{\ensuremath{\mathbb{#2}^{#3}_{#1}{}}}

  \providecommand{\FSpace}[3][]{\ensuremath{\ifx#2l \ell_{#3}^{#1}{}
      \else  #2_{#3}^{#1}{}\fi}} 
\providecommand{\SL}[1][2]{\ensuremath{\FSpace{SL}{#1}(\Space{R}{})}}
\providecommand{\norm}[2][\relax]{\left\|#2\right\|\ifx#1\relax\else_{#1}\fi}
\providecommand{\modulus}[2][\relax]{\left| #2 \right|\ifx#1\relax\else_{#1}\fi}
\providecommand{\scalar}[3][\relax]{\left\langle #2,#3 
        \right\rangle\ifx#1\relax\else_{#1}\fi}
\providecommand{\eqref}[1]{\textup{(\ref{#1})}}
\providecommand{\uir}[3][0]{\ifcase #1{\rho^{#2}_{#3}}%
\or {\breve{\rho}^{#2}_{#3}}%
\or {\tilde{\rho}^{#2}_{#3}}\fi}
\providecommand{\oper}[1]{\mathcal{#1}}
\providecommand{\rmi}{\mathrm{i}}
\providecommand{\MR}[1]{MR\#\href{http://www.ams.org/mathscinet-getitem?mr=#1}{#1}}
\providecommand{\Zbl}[1]{Zbl~\href{http://www.emis.de:80/cgi-bin/zmen/ZMATH/en/zmathf.html?first=1&maxdocs=3&type=html&an=#1&format=complete}{#1}}

\begin{document}
\title{Wavelets Beyond Admissibility}

\author
{\href{http://amsta.leeds.ac.uk/~kisilv/}{Vladimir V. Kisil}}

\address{%
School of Mathematics, 
University of Leeds, 
Leeds LS2\,9JT, 
UK\\
\href{mailto:kisilv@maths.leeds.ac.uk}{\texttt{kisilv@maths.leeds.ac.uk}}}

\begin{abstract} 
  The purpose of this paper is to articulate an
  observation that many interesting type of \emph{wavelets} (or
  \emph{coherent states}) arise from group representations which
  \emph{are not} square integrable or vacuum vectors which \emph{are
    not} admissible. 
\end{abstract}


\bodymatter

\section{Covariant Transform}
\label{sec:wavelet-transform}

A general group-theoretical
construction~\cite{Perelomov86,FeichGroech89a,Kisil98a,AliAntGaz00,Fuhr05a,ChristensenOlafsson09a}
of \emph{wavelets} (or \emph{coherent states}) starts from an
square integrable (\emph{s.i.})  representation
.  However, such a setup is
restrictive and is not necessary, in
fact. 

\begin{defn}
 Let \(\uir{}{}\) be a representation of
a group \(G\) in a space \(V\) and \(F\) be an operator from \(V\) to a space
\(U\). We define a \emph{covariant transform} 
\(\oper{W}\) from \(V\) to the space \(\FSpace{L}{}(G,U)\) of
\(U\)-valued functions on \(G\) by the formula:
\begin{equation}
  \label{eq:coheret-transf-gen}
  \oper{W}: v\mapsto \hat{v}(g) = F(\uir{}{}(g^{-1}) v), \qquad
  v\in V,\ g\in G.
\end{equation}
\end{defn}
\begin{rem}
  We do not require that operator \(F\) shall be linear.
\end{rem}
\begin{rem}
  \label{re:range-dim} 
  Usefulness of the covariant transform is in the reverse proportion
  to the dimensionality of the space \(U\). The covariant transform
  encodes properties of \(v\) in a function \(\oper{W}v\) on \(G\).
  For a low dimensional \(U\) this function can be ultimately
  investigated by means of harmonic analysis. Thus \(\dim U=1\) is the
  ideal case, however, it is unattainable sometimes, see
  Ex.~\ref{it:direct-funct} below.
\end{rem}
\begin{thm} 
  \label{pr:inter1} 
  The covariant transform \(\oper{W}\)~\eqref{eq:coheret-transf-gen}
  intertwines \(\uir{}{}\) and the left regular representation
  \(\Lambda\)    on \(\FSpace{L}{}(G,U)\):
  \begin{equation}\label{eq:left-reg-repr}
    \Lambda(g): f(h) \mapsto f(g^{-1}h).
  \end{equation}
\end{thm}
\begin{proof}
  We have a calculation similar to wavelet transform~\cite[Prop.~2.6]{Kisil98a}:
  \begin{displaymath}
    [\oper{W}( \uir{}{}(g) v)] (h) = F(\uir{}{}(h^{-1}) \uir{}{}(g) v ) 
    = [\oper{W}v] (g^{-1}h)
    = \Lambda(g) [\oper{W}v] (h).
  \end{displaymath}
\end{proof}
\begin{cor}\label{co:pi}
  The image space \(\oper{W}(V)\) is invariant under the
  left shifts on \(G\).
\end{cor}

\section{Examples of Covariant Transform}
\label{sec:exampl-covar-transf}

\begin{example}
  \label{ex:covariant}
  Let \(V\) be a Hilbert space with an inner product
  \(\scalar{\cdot}{\cdot}\) and \(\uir{}{}\) be a unitary
  representation. Let \(F: V \rightarrow \Space{C}{}\) be a
  functional \(v\mapsto \scalar{v}{v_0}\) defined by a vector
  \(v_0\in V\). Then the
  transformation~\eqref{eq:coheret-transf-gen} is the well-known
  expression for a \emph{wavelet transform}~\cite[(7.48)]{AliAntGaz00}
  (or \emph{representation coefficients}):
  \begin{equation}
    \label{eq:wavelet-transf}
    \oper{W}: v\mapsto \hat{v}(g) = \scalar{\uir{}{}(g^{-1})v}{v_0}  =
    \scalar{ v}{\uir{}{}(g)v_0}, \qquad
    v\in V,\ g\in G.
  \end{equation}
  The family of vectors \(v_g=\uir{}{}(g)v_0\) is called
  \emph{wavelets} or \emph{coherent states}. In this case we obtain
  scalar valued functions on \(G\), thus the fundamental r\^ole of
  this example is explained in Rem.~\ref{re:range-dim}.

  This scheme is typically carried out for a s.i. representation
  \(\uir{}{}\) and \(v_0\) being an admissible
  vector\cite{Perelomov86,FeichGroech89a,AliAntGaz00,Fuhr05a,ChristensenOlafsson09a}.
  In this case the wavelet (covariant) transform is a map into the s.i.
  functions~\cite{DufloMoore} with respect to the left Haar
  measure.
\end{example}
However s.i. representations and admissible vectors does not cover all
interesting cases.
\begin{example}
  \label{ex:ax+b}
  Let \(G\) be the ``\(ax+b\)'' (or \emph{affine})
  group~\cite[\S~8.2]{AliAntGaz00}: the set of points \((a,b)\),
  \(a\in \Space[+]{R}{}\), \(b\in \Space{R}{}\) in the upper
  half-plane with the group law:
   \begin{equation}
     (a, b) * (a', b') = (aa', ab'+b)
   \end{equation}
   and left invariant measure \(a^{-2}\,da\,db\).  Its isometric
   representation on \(V=\FSpace{L}{p}(\Space{R}{})\) is given by the
   formula:
  \begin{equation}\label{eq:ax+b-repr-1}
     [\uir{}{p}(a,b)\, f](x)= a^{\frac{1}{p}}f\left(ax+b\right).
  \end{equation}
  We consider the operators \(F_{\pm}:\FSpace{L}{2}(\Space{R}{})
  \rightarrow \Space{C}{}\) defined by:
  \begin{equation}
    \label{eq:cauchy}
    F_{\pm}(f)=\frac{1}{2\pi i}\int_{\Space{R}{}} \frac{f(t)\,dt}{t\mp \rmi}.
  \end{equation}
  Then the covariant transform~\eqref{eq:coheret-transf-gen} is the
  Cauchy integral from \(\FSpace{L}{2}(\Space{R}{})\) to the Hardy
  space in the upper/lower half-plane
  \(\FSpace{H}{2}(\Space[\pm]{R}{2})\).  Although the
  representation~\eqref{eq:ax+b-repr-1} is s.i. for \(p=2\), the
  function \(\frac{1}{t\pm \rmi}\) is not an admissible vacuum vector.
  Thus the complex analysis become decoupled from the traditional
  wavelets theory. As a result the application of wavelet theory shall
  relay on an extraneous mother wavelets~\cite{Hutnik09a}.

  However many important objects in complex analysis are generated by
  inadmissible mother wavelets like~\eqref{eq:cauchy}. For example, if
  \(F:\FSpace{L}{2}(\Space{R}{}) \rightarrow \Space{C}{}\) is defined
  by \(F: f \mapsto F_+ f + F_-f\) then the covariant
  transform~\eqref{eq:coheret-transf-gen} is simply the \emph{Poisson
    integral}.  If \(F:\FSpace{L}{2}(\Space{R}{}) \rightarrow
  \Space{C}{2}\) is defined by \(F: f \mapsto( F_+ f, F_-f)\) then the
  covariant transform~\eqref{eq:coheret-transf-gen} represents a
  function on the real line as a jump between functions analytic in
  the upper and the lower half-planes. This makes a decomposition of
  \(\FSpace{L}{2}(\Space{R}{})\) into irreducible components of the
  representation~\eqref{eq:ax+b-repr-1}.  Another interesting but
  non-admissible vector is the Gaussian \(e^{-x^2}\).
\end{example}
\begin{example}
  \label{ex:sl2}
  For the group \(G=\SL\)~\cite{Lang85} let us consider the unitary
  representation \(\uir{}{}\) on the space of s.i. function
  \(\FSpace{L}{2}(\Space[+]{R}{2})\) on the upper half-plane through
  the M\"obius transformations:
   \begin{displaymath}
     \uir{}{}(g): f(z) \mapsto \frac{1}{(c z + d)^2}\,
     f\left(\frac{a z+ b }{c z +d}\right), \qquad g^{-1}=\
     \begin{pmatrix}
       a & b \\ c & d 
     \end{pmatrix}.
   \end{displaymath}
   Let \(F_i\) be the functional
   \(\FSpace{L}{2}(\Space[+]{R}{2})\rightarrow \Space{C}{}\) of
   pairing with the lowest/highest \(i\)-weight vector in the
   corresponding irreducible component of the discrete
   series~\cite[Ch.~VI]{Lang85}. Then we can build an operator \(F\)
   from various \(F_i\) similarly to the previous example, e.g.  this
   generalises the representation of an s.i. function as a sum of
   analytic ones from different irreducible subspaces.

   Covariant transform is also meaningful for principal and
   complementary series of representations of the group
   \(\SL\)\cite{Kisil97c}, which are not s.i. 
\end{example}
\begin{example}
  \label{it:direct-funct}
  A straightforward generalisation of Ex.\ref{ex:covariant} is
  obtained if \(V\) is a Banach space and \(F: V \rightarrow
  \Space{C}{}\) is an element of \(V^*\). Then the
  covariant transform coincides with the construction of wavelets in
  Banach spaces~\cite{Kisil98a}.

  The next stage of generalisation is achieved if \(V\) is a
  Banach space and \(F: V \rightarrow \Space{C}{n}\) be a linear
  operator. Then the corresponding covariant transform is a map
  \(\oper{W}: V \rightarrow \FSpace{L}{}(G,\Space{C}{n})\). This is
  closely related to M.G.~Krein's works on \emph{directing
    functionals}~\cite{Krein48a}, see also \emph{multiresolution
    wavelet analysis}~\cite{BratJorg97a}, Clifford-valued
  Bargmann spaces~\cite{CnopsKisil97a} and~\cite[Thm.~7.3.1]{AliAntGaz00}.
\end{example}
\begin{example}
  A step in a different direction is a consideration of
  non-linear operators. Take again the ``\(ax+b\)'' group and its
  representation~\eqref{eq:ax+b-repr-1}. 
  We define \(F\) to be a homogeneous but non-linear functional
  \(V\rightarrow \Space[+]{R}{}\):
  \begin{displaymath}
    F (f) = \frac{1}{2}\int\limits_{-1}^1 \modulus{f(x)}\,dx.
  \end{displaymath}
  The covariant transform~\eqref{eq:coheret-transf-gen} becomes:
  \begin{eqnarray*}
    [\oper{W}_p f](a,b) 
    = \frac{1}{2}\int\limits_{-1}^1
    \modulus{a^{\frac{1}{p}}f\left(ax+b\right)}\,dx
    = a^{\frac{1}{p}}\frac{1}{2a}\int\limits^{b+a}_{b-a}
    \modulus{f\left(x\right)}\,dx.
  \end{eqnarray*}
  Obviously \(M_f(b)=\max_{a}[\oper{W}_{\infty}f](a,b)\) coincides
  with the Hardy \emph{maximal function}, which contains important
  information on the original function \(f\). However, the full
  covariant transform is even more detailed. For example,
  \(\norm{f}=\max_b[\oper{W}_{\infty}f](\frac{1}{2},b)\) is the
  shift invariant norm~\cite{Johansson08a}. 

  From the Cor.~\ref{co:pi} we deduce that the operator \(M: f\mapsto
  M_f\) intertwines \(\uir{}{p}\) with itself \(\uir{}{p}M=M
  \uir{}{p}\).
\end{example}
\begin{example}
  Let \(V=\FSpace{L}{c}(\Space{R}{2})\) be the space of
  compactly supported bounded functions on the plane. We take \(F\)
  be the linear operator \(V\rightarrow \Space{C}{}\) of integration
  over the real line:
  \begin{displaymath}
    F: f(x,y)\mapsto F(f)=\int_{\Space{R}{}}f(x,0)\,dx.
  \end{displaymath}
  Let \(G\) be the group of Euclidean motions of the plane
  represented by \(\uir{}{}\) on \(V\) by a change of variables. Then
  the wavelet transform \(F(\uir{}{}(g)f)\) is the \emph{Radon
    transform}.
\end{example}
\begin{example}
  Let a representation \(\uir{}{}\) of a group \(G\) act on a
  space \(X\). Then there is an associated representation
  \(\uir{}{B}\) of \(G\) on a space \(V=\FSpace{B}{}(X,Y)\) of
  linear operators \(X\rightarrow Y\) defined by the identity:
  \begin{displaymath}
    (\uir{}{B}(g) A)x=A(\uir{}{}(g)x), \qquad x\in X,\ g\in G,\ A \in \FSpace{B}{}(X,Y).
  \end{displaymath}
  Following the Remark~\ref{re:range-dim} we take \(F\) to be a
  functional \(V\rightarrow\Space{C}{}\), for example \(F\) can be
  defined from a pair \(x\in X\), \(l\in Y^*\)  by the expression
  \(F: A\mapsto \scalar{Ax}{l}\). Then the covariant
  transform:
  \begin{displaymath}
    \oper{W}: A \mapsto \hat{A}(g)=F(\uir{}{B}(g) A), \qquad 
  \end{displaymath}
  this is an example of \emph{covariant calculus}~\cite{Kisil98a,Kisil02a}.
\end{example}
\begin{example}
  A modification of the previous construction is obtained if we
  have two groups \(G_1\) and \(G_2\) represented by \(\uir{}{1}\)
  and \(\uir{}{2}\) on \(X\) and \(Y^*\) respectively. Then we have a covariant
  transform \(\FSpace{B}{}(X,Y)\rightarrow \FSpace{L}{}(G_1\times
  G_2, \Space{C}{})\) defined by the formula:
  \begin{displaymath}
    \oper{W}: A \mapsto \hat{A}(g_1,g_2)=\scalar{A\uir{}{1}(g_1)x}{\uir{}{2}(g_2)l}.
  \end{displaymath}
  This generalises \emph{Berezin functional calculi}~\cite{Kisil98a}.
\end{example}
\begin{example}
  Let us restrict the previous example to the case when \(X=Y\) is a
  Hilbert space, \(\uir{}{1}{}=\uir{}{2}{}=\uir{}{}\) and \(x=l\) with
  \(\norm{x}=1\). Than the range of the covariant transform:
  \begin{displaymath}
    \oper{W}: A \mapsto \hat{A}(g)=\scalar{A\uir{}{}(g)x}{\uir{}{}(g)x}
  \end{displaymath}
  is a subset of the \emph{numerical range} of the operator \(A\). 
\end{example}
\begin{example}
  The group \(\SL\) consists of \(2\times 2\) matrices
  of the form \(\begin{pmatrix}
    \alpha&\beta\\\bar{\beta}&\bar{\alpha}
  \end{pmatrix}\) with the unit
  determinant~\cite[\S~IX.1]{Lang85}. Let \(A\) be an operator
  with the spectral radius less than \(1\). Then the associated
  M\"obius transformation
  \begin{displaymath}
    g: A \mapsto g\cdot A =  \frac{\alpha A+\beta
      I}{\bar{\beta}A+\bar{\alpha}I}, \qquad \text{where} \quad
    g^{-1}=
    \begin{pmatrix}
      \alpha&\beta\\\bar{\beta}&\bar{\alpha}
    \end{pmatrix}\in \SL,\ 
  \end{displaymath}
  produces a well-defined operator with the spectral radius less than
  \(1\) as well.  Thus we have a representation of \(\SL\). A choise
  of an operator \(F\) will define the corresponding covariant
  transform. In this way we obtain generalisations of
  \emph{Riesz--Dunford functional calculus}~\cite{Kisil02a}.
\end{example}

\section{Inverse Covariant Transform}
\label{sec:invar-funct-groups}
An object invariant under the left action
\(\Lambda\)~\eqref{eq:left-reg-repr} is called \emph{left invariant}.
For example, 
let \(L\) and \(L'\) be two left invariant spaces of functions on
\(G\).  We say that a pairing \(\scalar{\cdot}{\cdot}: L\times
L' \rightarrow \Space{C}{}\) is \emph{left invariant} if
\begin{equation}
  \scalar{\Lambda(g)f}{\Lambda(g) f'}= \scalar{f}{f'}, \quad \textrm{ for all }
  \quad f\in L,\  f'\in L'.
\end{equation}
\begin{rem}
  \begin{enumerate}
  \item We do not require the pairing to be linear in general.
  \item If the pairing is invariant on space \(L\times L'\) it is not
    necessarily invariant (or even defined) on the whole
    \(\FSpace{C}{}(G)\times \FSpace{C}{}(G)\).
  \item In a more general setting we shall study an invariant pairing
    on a homogeneous spaces instead of the group. However due to length
    constraints we cannot consider it here beyond the Example~\ref{ex:hs-pairing}.
  \item An invariant pairing on \(G\) can be obtained from an invariant
    functional \(l\) by the formula \(\scalar{f_1}{f_2}=l(f_1\bar{f}_2)\).
  \end{enumerate}
\end{rem}

For a representation \(\uir{}{}\) of \(G\) in \(V\) and \(v_0\in V\)
we fix a function \(w(g)=\uir{}{}(g)v_0\). We assume that the pairing
can be extended in its second component to this \(V\)-valued
functions, say, in the weak sense.
\begin{defn}
  \label{de:admissible}
  Let \(\scalar{\cdot}{\cdot}\) be a left invariant pairing on
  \(L\times L'\) as above, let \(\uir{}{}\) be a representation of
  \(G\) in a space \(V\), we define the function
  \(w(g)=\uir{}{}(g)v_0\) for \(v_0\in V\). The \emph{inverse
    covariant transform} \(\oper{M}\) is a map \(L \rightarrow V\)
  defined by the pairing:
  \begin{equation}
    \label{eq:inv-cov-trans}
    \oper{M}: f \mapsto \scalar{f}{w}, \qquad \text{
      where } f\in L. 
  \end{equation}
\end{defn}

 \begin{example}
   Let \(G\) be a group with a unitary s.i. representation \(\rho\).
   An invariant pairing of two s.i. functions is obviously done by the
   integration over the Haar measure:
   \begin{displaymath}
     \scalar{f_1}{f_2}=\int_G f_1(g)\bar{f}_2(g)\,dg.
   \end{displaymath}
   
   For an admissible vector \(v_0\)~\cite{DufloMoore},
   \cite[Chap.~8]{AliAntGaz00} the inverse covariant transform is
   known in this setup as \emph{reconstruction formula}.
 \end{example}
 \begin{example}
   \label{ex:hs-pairing}
   Let \(\rho\) be a s.i. representation of \(G\) modulo a subgroup
   \(H\subset G\) and let \(X=G/H\) be the corresponding homogeneous
   space with a quasi-invariant measure \(dx\).  Then integration over
   \(dx\) with an appropriate weight produces an invariant pairing.
   The inverse covariant transform is a more general
   version~\cite[(7.52)]{AliAntGaz00} of the \emph{reconstruction
     formula} mentioned in the previous example.
 \end{example}
 

 Let \(\rho\) be not a s.i. representation (even modulo a subgroup) or
 let \(v_0\) be inadmissible vector of a s.i.  representation
 \(\rho\). An invariant pairing in this case is not associated with an
 integration over any non singular invariant measure on \(G\). In this
 case we have a \emph{Hardy pairing}. The following example explains
 the name.
\begin{example}
  Let \(G\) be the ``\(ax+b\)'' group and its representation
  \(\uir{}{}\)~\eqref{eq:ax+b-repr-1} from Ex.~\ref{ex:ax+b}. An
  invariant pairing on \(G\), which is not generated by the Haar
  measure \(a^{-2}da\,db\), is:
  \begin{equation}
    \label{eq:hardy-pairing}
    \scalar{f_1}{f_2}=
    \lim_{a\rightarrow 0}\int\limits_{-\infty}^{\infty}
    f_1(a,b)\,\bar{f}_2(a,b)\,db.
  \end{equation}
  For this pairing we can consider functions \(\frac{1}{2\pi i
    (x+i)}\) or \(e^{-x^2}\), which are not admissible vectors in the
  sense of s.i. representations. Then the inverse covariant transform
  provides an \emph{integral resolutions} of the identity.
\end{example}
Similar pairings can be defined for other semi-direct products of two
groups. We can also extend a Hardy pairing to a group, which has a
subgroup with such a pairing.
\begin{example}
  Let \(G\) be the group \(\SL\) from the Ex.~\ref{ex:sl2}. Then
  the ``\(ax+b\)'' group is a subgroup of \(\SL\), moreover we can
  parametrise \(\SL\) by triples \((a,b,\theta)\),
  \(\theta\in(-\pi,\pi]\) with the respective Haar
  measure~\cite[III.1(3)]{Lang85}. Then the Hardy
  pairing
  \begin{equation}
    \label{eq:hardy-pairing}
    \scalar{f_1}{f_2}= \lim_{a\rightarrow 0}\int\limits_{-\infty}^{\infty}
    f_1(a,b,\theta)\,\bar{f}_2(a,b,\theta)\,db\,d\theta.
  \end{equation}
  is invariant on \(\SL\) as well.  The corresponding inverse
  covariant transform provides even a finer resolution of the identity
  which is invariant under conformal mappings of the Lobachevsky
  half-plane.
\end{example}

A further study of covariant transform and its inverse shall be
continued elsewhere.

\small

\begin{thebibliography}{10}

\bibitem{Perelomov86}
A.~Perelomov, {\em Generalized coherent states and their applications}
(Springer-Verlag,
  Berlin, 1986).

\bibitem{FeichGroech89a}
{Feichtinger, Hans G. and Groechenig, K.H.}, {\em J. Funct. Anal.} {\bf 86},
  307 (1989)
.

\bibitem{Kisil98a}
V.~V. Kisil, {\em Acta Appl. Math.} {\bf 59}, 79 (1999),
\arXiv{math/9807141}
.

\bibitem{AliAntGaz00}
S.~T. Ali, J.-P. Antoine and J.-P. Gazeau, {\em Coherent States, Wavelets and
  Their Generalizations} (Springer-Verlag, New York, 2000).

\bibitem{Fuhr05a}
H.~F{\"u}hr, {\em Abstract Harmonic Analysis of Continuous Wavelet Transforms},
  Lecture Notes in Mathematics, Vol.~1863 (Springer-Verlag, Berlin, 2005).

\bibitem{ChristensenOlafsson09a}
J.~G. Christensen and G.~{\'O}lafsson, {\em Acta Appl. Math.} {\bf 107}, 25
  (2009).

\bibitem{DufloMoore}
M.~Duflo and C.~C. Moore, {\em J. Functional Analysis} {\bf 21}, 209 (1976).

\bibitem{Hutnik09a}
O.~Hutn{\'{\i}}k, {\em Integral Equations Operator Theory} {\bf 63}, 29 (2009).

\bibitem{Kisil97c}
V.~V. Kisil, {\em Complex Variables Theory Appl.} {\bf 40}, 93 (1999),
  \arXiv{funct-an/9712003}.

\bibitem{Krein48a}
M.~G. Kre{\u\i}n, {\em Akad. Nauk Ukrain. RSR. Zbirnik Prac\cprime\ Inst. Mat.}
  {\bf 1948}, 83 (1948), \MR{14:56c}, reprinted in~\cite{KreinII}.

\bibitem{BratJorg97a}
O.~Bratteli and P.~E.~T. Jorgensen, {\em Integral Equations Operator Theory}
  {\bf 28}, 382 (1997), \arXiv{funct-an/9612003}.

\bibitem{CnopsKisil97a}
J.~Cnops and V.~V. Kisil, {\em Math. Methods Appl. Sci.} {\bf 22}, 353 (1999),
  \arXiv{math/9806150}. \Zbl{1005.22003}.

\bibitem{Johansson08a}
A.~Johansson, {\em Systems Control Lett.} {\bf 57}, 105 (2008).

\bibitem{Lang85}
S.~Lang, {\em {${\rm SL}\sb 2({\bf R})$}} (Springer-Verlag, New York, 1985).


\bibitem{Kisil02a}
V.~V. Kisil, Spectrum as the support of functional calculus, in {\em Functional
  analysis and its applications\/}, North-Holland
Math. Stud. Vol.~197, pp.~133--141,
  (Elsevier, Amsterdam, 2004).
\newblock \arXiv{math.FA/0208249}.

\bibitem{KreinII}
M.~G. Kre{\u\i}n, {\em {\cyr {I}zbrannye Trudy}. {II}} (Akad. Nauk Ukrainy
  Inst. Mat., Kiev, 1997).
\newblock \MR{96m:01030}.

\end{thebibliography}
\newcommand{\noopsort}[1]{} \newcommand{\printfirst}[2]{#1}
  \newcommand{\singleletter}[1]{#1} \newcommand{\switchargs}[2]{#2#1}
  \newcommand{\irm}{\textup{I}} \newcommand{\iirm}{\textup{II}}
  \newcommand{\vrm}{\textup{V}} \providecommand{\cprime}{'}
  \providecommand{\eprint}[2]{\texttt{#2}}
  \providecommand{\myeprint}[2]{\texttt{#2}}
  \providecommand{\arXiv}[1]{\myeprint{http://arXiv.org/abs/#1}{arXiv:#1}}
  \providecommand{\CPP}{\texttt{C++}} \providecommand{\NoWEB}{\texttt{noweb}}
  \providecommand{\MetaPost}{\texttt{Meta}\-\texttt{Post}}
  \providecommand{\GiNaC}{\textsf{GiNaC}}
  \providecommand{\pyGiNaC}{\textsf{pyGiNaC}}
  \providecommand{\Asymptote}{\texttt{Asymptote}}

\end{document}